\documentclass{amsart}
\usepackage{amsmath,amssymb,amsthm,verbatim}
\usepackage[pdftex,bookmarks=true]{hyperref}
\newtheorem{theorem}{Theorem}[section] % 1st argument is your name for it
\newtheorem{lemma}[theorem]{Lemma}     % 2nd argument is what is printed

\title% end with percent
 {Reducing the number of prime factors of long $\kappa$-tuples} % This is the full title of the paper
\author{C. S. Franze}

\begin{document}
\maketitle

\begin{abstract}
We prove that there are infinitely many integers $n$ such that the total number of prime factors of $(n+h_{1})\ldots(n+h_{\kappa})$ is at most $\frac{1}{2}\kappa\log\kappa+O(\kappa)$, provided $\kappa$ is sufficiently large.
\end{abstract}

%\part{Use this type of header for very long papers only}
% use lowercase except for proper names

\section{Introduction} % use lowercase except for proper names
\noindent In this paper, we present a weighted sieve method and apply it to obtain an improvement in $\kappa$-tuples with few prime factors. More specifically, we show that an admissible tuple of length $\kappa$ is infinitely often a product of at most $r_{\kappa}$ prime factors, where $r_{\kappa}\sim\frac{1}{2}\kappa\log\kappa$, provided $\kappa$ is sufficiently large. This result is stated more carefully in Theorem \ref{Th:1} below. Until now, the best known results had $r_{\kappa}\sim\kappa\log\kappa$. Of course, if the $\kappa$-tuples conjecture is true, then we may take $r_{\kappa}=\kappa$. We obtain our results by considering sums of the form
\begin{equation} \label{example1}
\sum_{\substack{n\in \mathcal{A}\\(n,P(y))=1}}
\left( \sum_{\substack{d|n\\ d|P(z)}} a_d \right)
\left( \sum_{\substack{\nu|n \\ \nu|P(z')}} \lambda_\nu\right)^2.
\end{equation}
We will suppose that $z'\le z$, and in fact, for this application at least, good results are obtained when $z'$ is roughly $z^{1/2\kappa}$. The choice for $a_d$ will be inspired by the Richert weights, and the choice of $\lambda_\nu$ will be motivated by Selberg's upper bound sieve.\\
\indent The first published application of a sieve construction such as \eqref{example1} to almost-primes generated by a polynomial
was given by R. Miech \cite{M64}, who used Kuhn's weights and a classical choice for the $\lambda_\nu$'s.
Miech made use of unpublished notes of I. Reiner, P. Bateman, and L. Rubel of Selberg's lectures given at
the Institute for Advanced Study in 1948, 1950, and 1958. Halberstam and Richert \cite[Section 10.4]{HR74} gave an account that employed Richert's weights. In this paper, we also use Richert's weights, but employ a non-trivial range for the sieve function $j_\kappa$. We also allow for a non-classical choice of the $\lambda_\nu$'s, a helpful generalization for future work on reducing the number of prime factors of \emph{short} $\kappa$-tuples.\\
\indent Our result in Theorem \ref{Th:1} represents an improvement over the work of Miech \cite{M64}, who showed that one could take $r_{\kappa}\sim\kappa\log\kappa$. More recently, Ramar\'e \cite{Ramare10} was able to show that there are many $\kappa$-tuples with \textit{exactly} $\kappa\log\kappa+O(\kappa\sqrt{\log\kappa})$ prime factors, suggesting that obtaining $\kappa$-tuples with $(1-\delta)\kappa\log\kappa$ prime factors, for any $\delta>0$, would be a challenge. In another direction, Ho and Tsang \cite{HT06}, following the work of Heath-Brown, found that one could take $r_{\kappa}\sim1.44\ \kappa\log\kappa$. However, their emphasis was different in that they sought to minimize the number of prime factors occurring in \textit{each} of the terms of $L(n)$.\\

\begin{theorem}\label{Th:1}
Suppose $L(n)$ is a product of $\kappa$ linear forms,
\begin{equation}\label{Ldef}
L(n)=\prod_{i=1}^{\kappa}\left(a_{i}n+b_{i}\right),
\end{equation}
$(a_{i},b_{i})=1$, and $\Delta_{L}\neq0$, where $\Delta_{L}$ is the discriminant of $L(n)$, defined as
\begin{equation}\label{E:disc}
\Delta_{L}=\prod_{i=1}^{\kappa}a_{i}\prod_{1\le t<s\le\kappa}(a_{t}b_{s}-a_{s}b_{t}).
\end{equation}
Define $\rho(p)$ to be the number of solutions to $L(n)\equiv0\mod p$, and suppose that $\rho(p)<p$ for all primes $p\le\kappa$. Then, for all sufficiently large $\kappa$ and $x$, we have
\begin{equation*}
\left|\left\{n\le x: \Omega(L(n))\le r_{\kappa}\right\}\right|\gg \frac{x}{\log^{\kappa}x},
\end{equation*}
for any $r_{\kappa}$ satisfying the inequality
\begin{equation*}
r_{\kappa}>\frac{1}{2}\kappa\log\kappa+\left(1+\frac{\gamma}{2}+\log 4\right)\kappa+\frac{13}{18}\sqrt{\frac{\kappa}{\pi}}+O\left(\log\kappa\right).
\end{equation*}
\end{theorem}
\indent More generally, one could use our construction to consider a polynomial $H(n)$ which is the product of $\kappa$ irreducible polynomials, each of degree $h$,
\begin{equation*}
H(n)=\prod_{i=1}^{\kappa}h_{i}(n).
\end{equation*}
In this case, one can show that for sufficiently large $\kappa$,
\begin{equation*}
r_{\kappa}\sim \kappa h+\frac{1}{2}\kappa\log\kappa+\kappa\log 2h
\end{equation*}
is an admissible choice for $r_{\kappa}$. When $h$ is small compared to $\kappa$, this bound is superior to that given by other constructions, provided $\kappa$ is taken sufficiently large. On the other hand, if $h$ is large compared to $\kappa$, the $\log 2h$ term gets out of control. This phenomenon is noted in the case $\kappa=1$ in Halberstam and Richert \cite[Section 10.5]{HR74}. For instance, using our construction it is possible to show that $r_{1}=h+1+\log 2h$ is an admissible choice. However, other constructions take full advantage of the linear sieve to obtain $r_{1}=h+1$. For this reason, we focus on the case when $h=1$ and $\kappa$ is large.
\section{Preliminaries}

In this section, we wish to provide the reader with the relevant framework associated to this particular sieve problem. In the general setting, one starts with a finite sequence of integers, say $\mathcal{A}$, and a set of primes, say $\mathcal{P}$. In our case, we will take
\begin{align*}
\mathcal{A}=\left\{L(n):n\le x\right\},
\end{align*}
where $L(n)$ is a product of $\kappa$ linear forms, as in \eqref{Ldef}, and $\mathcal{P}$ to be the set of all primes $p$ less than $z$. For future reference, we define
\begin{equation*}
P(z)=\prod_{p<z}p.
\end{equation*}

Following usual notation, we let $\mathcal{A}_{d}$ be the elements of $\mathcal{A}$ that are divisible by $d$. The first step is to understand $\left|\mathcal{A}_{d}\right|$, the size of $\mathcal{A}_{d}$. Actually, we will need only to understand $\left|\mathcal{A}_{d}\right|$ for squarefree $d$. One typical sieve assumption in this direction is that there exists a multiplicative function, say $f$, such that
\begin{equation}\label{sieveassumption}
|\mathcal{A}_d|=\frac{X}{f(d)}+\mathcal{R}_{d},
\end{equation}
where the $\mathcal{R}_{d}$ are small, at least on average. In our example we have
\begin{equation*}
\left|\mathcal{A}_{d}\right|=\sum_{\substack{n\le x\\L(n)\equiv0(d)}}1=\rho(d)\left(\frac{x}{d}+\theta\right),
\end{equation*}
where $\left|\theta\right|\le 1$, and $\rho(d)$ is the number of solutions to $L(n)\equiv0$ mod $d$. This implies that $\left|\mathcal{R}_{d}\right|\le\rho(d)$. Furthermore, an application of the Chinese Remainder Theorem shows that $\rho(d)$ is multiplicative (recall that $d$ is squarefree), so that
\begin{align*}
\frac{\rho(d)}{d}=\frac{1}{f(d)}.
\end{align*}
The condition that $(a_{i},b_{i})=1$ appearing in Theorem \ref{Th:1} guarantees that $\rho(d)\neq0$, and hence $f(d)$ is well-defined. This is enough to show that the assumption in \eqref{sieveassumption} is valid, and that $X=x$.\\
\indent Next, we will outline how sums such as \eqref{example1} are dealt with under the minimal assumptions above. In particular, there is a clever choice for the $\lambda_{\nu}$ that allow these sums to be decomposed into a main term and an error term. To begin with, let $f'=f*\mu$, and $\lambda_{\nu}$ be an arbitrary sequence of real numbers with the property that $\lambda_\nu=0$ if $\nu\nmid P(z')$ or if $\nu > \xi $. Assume that $\lambda_{1}\neq0$. Now, define a new sequence $\zeta_r$ by the relation
\begin{equation*}
\frac{\mu(r)\zeta_r}{f'(r)} = \sum_{d<\frac{\xi}{r}} \frac{\lambda_{dr}}{f(dr)}.
\end{equation*}
By M\"obius inversion, we also have
\begin{equation*}
\frac{\mu(d)\lambda_d}{f(d)} = \sum_{r<\frac{\xi}{d}} \frac{\zeta_{dr}}{f'(dr)}.
\end{equation*}
Having made these assumptions of $\lambda_{\nu}$, we have the identity
\begin{align}\label{E:Sieve}
\sum_{n\in\mathcal{A}}\left(\sum_{\substack{d|n\\d\mid P(z)}}a_d\right)\left(\sum_{\substack{\nu|n\\\nu\mid P(z')}}\lambda_\nu \right)^2=X \mathfrak{S}_{\mathcal{A}}+\mathfrak{E}_{\mathcal{A}},
\end{align}
where
\begin{align}\label{E:Main}
\mathfrak{S}_{\mathcal{A}}=\sum_{\substack{m<\xi\\m\mid P(z')}} \sum_{\substack{d\mid P(z)\\(d,m)=1}}\frac{\mu^2(m)}{f'(m)}\frac{a_d}{f(d)}\left( \sum_{r|d} \mu(r) \zeta_{rm} \right)^{2},
\end{align}
and
\begin{align}\label{E:Error}
\mathfrak{E}_{\mathcal{A}}=\sum_{\substack{d\mid P(z)\\\nu_{1},\nu_{2}\mid P(z')}}a_{d}\lambda_{\nu_{1}}\lambda_{\nu_{2}}\mathcal{R}_{\left[d,\nu_{1},\nu_{2}\right]}.
\end{align}
This identity is the starting point of Selberg's lower bound sieve method, and, in the case when $z=z'$, has appeared in Selberg \cite[see Section 7 on p.82]{Selberg}, Bombieri \cite[see Theorem 18 on p.65]{Bombieri}, Cojocaru and Murty \cite[see Theorem 10.11 on p.178]{CoM}, Greaves \cite[see Lemma 1 on p.286]{Greaves}, and others. A trivial modification allows for the case when $z\neq z'$. If $\mathfrak{S}_{\mathcal{A}}$ remains positive even as $X\rightarrow\infty$, then the sieve will be successful at achieving a positive lower bound, provided the error term $\mathfrak{E}_{\mathcal{A}}$ is negligible.\\
\indent In analyzing $\mathfrak{S}_{\mathcal{A}}$, we will encounter the well-known sieve quantity
\begin{align*}
V(z')=\prod_{p<z'}\left(1-\frac{\rho(p)}{p}\right).
\end{align*}
Recall that in Theorem \ref{Th:1} we assume that $\rho(p)<p$ for all primes $p$. Therefore, we have that $V(z')\neq0$. In fact, one can easily verify that
\begin{equation}\label{B:Vbound}
\frac{1}{V(z')}\ll\log^{\kappa}z',
\end{equation}
and
\begin{equation}\label{E:HSum}
\sum_{p<s}\frac{\rho(p)\log p}{p}=\kappa\log s+O\left(1\right),
\end{equation}
since $\rho(p)=\kappa$ for most primes $p$.

\section{The Richert weights}
The Richert weights are defined by
\begin{equation}\label{W:richertweights}
	a_d =
	\begin{cases}
	   \phantom{+} b & \text{if $d=1$,}\\
	    -b & \text{if $d$ is prime and $d<y$,} \\
	    -\frac{\log(\frac{z}{d})}{\log z} & \text{if $d$ is prime and $y\le d<z$,} \\
	    \phantom{+}0 & \text{otherwise}.
	 \end{cases}
\end{equation}
The weight attached to the primes $p<y$ is a device that allows us to remove the condition that $(n,P(y))=1$ appearing in \eqref{example1}. Ultimately we will end up taking $y$ to be a very small power of $x$. Ignoring the contribution from these primes, and furthermore removing the $n\in\mathcal{A}$ that are divisible by the square of a prime $p\in\left[y,z\right)$, the Richert weights allow us to bound $\Omega(n)$ using $b$. For, in this case, if
\begin{equation*}
\sum_{\substack{d\mid n\\d\mid P(z)}}a_{d}>0,
\end{equation*}
then
\begin{equation*}
0<b-\sum_{\substack{p\mid n\\y\le p<z}}\left(1-\frac{\log p}{\log z}\right)\le b-\Omega(n)+\frac{\log\left|n\right|}{\log z}.
\end{equation*}

The goal of this section is to use these weights to prove
\begin{lemma}\label{L:weightedsieve}
Suppose $L(n)$ is subject to the hypotheses of Theorem \ref{Th:1}. Then, for all sufficiently large $x$, and any $\varepsilon>0$, we have
\begin{equation*}
\left|\left\{n\le x: \Omega(L(n))\le r_{\kappa}\right\}\right|\gg \frac{x}{\log^{\kappa}x}\Big(\mathfrak{S}_{\mathcal{A}}V(z')+o(1)\Big),
\end{equation*}
provided $r_{\kappa}>U\kappa-1+O(\varepsilon)$, where $z=x^{\frac{1}{U}}$, $z'=\xi^{\frac{1}{u}}$, and $z\xi^{2}=x^{1-\delta}$, for any $\delta>0$.
\end{lemma}
The trick to using Lemma \ref{L:weightedsieve} is to choose $U$ so that $\mathfrak{S}_{\mathcal{A}}V(z')>0$ \textit{and} $r_{\kappa}$ remains small. An innovation employed by Halberstam and Richert \cite{HR74} allows for such a choice of $U$ and will be discussed in Section \ref{innovation}.\\
\indent In preparation for the proof of Lemma \ref{L:weightedsieve}, some comments concerning the $\lambda_{\nu}$'s appearing in $\mathfrak{S}_{\mathcal{A}}$ are in order. First, recall that
\begin{equation*}
\frac{\mu(d)\lambda_d}{f(d)} = \sum_{r<\frac{\xi}{d}} \frac{\zeta_{dr}}{f'(dr)}.
\end{equation*}
The $\zeta_{r}$ will be chosen as
\begin{equation}\label{E:polyP}
\zeta_{r}=P^{*}\left(\frac{\log\xi/r}{\log z'}\right)
\end{equation}
where $P^{*}(w)$ is a polynomial that is positive for $0\le w\le u$. Therefore,
\begin{equation}\label{B:firstlambdabound}
\lambda_{1}=\sum_{r<\xi}\frac{\zeta_{r}}{f'(r)}\le\sup_{0\le w\le u}P(w)\sum_{\substack{r<\xi\\r\mid P(z')}}\frac{1}{f'(r)}\ll\sum_{r\mid P(z')}\frac{1}{f'(r)}=\frac{1}{V(z')}.
\end{equation}
In the case when $\zeta_{r}=1$, the $\lambda_{v}$ are well understood. We will refer to this choice of $\lambda_{v}$ as $\widetilde{\lambda_{v}}$. It is known, for example, that $\left|\widetilde{\lambda_{v}}\right|\le\left|\widetilde{\lambda_{1}}\right|$. A proof of this fact can be found in Halberstam and Diamond \cite[Section 2.2]{DHG08}. Since
\begin{equation*}
\left|\lambda_{\nu}\right|\le\sup_{0\le w\le u}P(w)\widetilde{\lambda_{1}},
\end{equation*}
and
\begin{equation*}
\lambda_{1}=\sum_{\substack{r<\xi\\r\mid P(z')}}\frac{\mu^2(r)}{f'(r)}P\left(\frac{\log\xi/r}{\log z'}\right)\ge\inf_{0\le w\le u}P(w)\widetilde{\lambda_{1}},
\end{equation*}
it is clear that
\begin{equation}\label{B:lambdabound}
\frac{\left|\lambda_{v}\right|}{\left|\lambda_{1}\right|}\le\displaystyle\frac{\displaystyle\sup_{0\le w\le u}P(w)}{\displaystyle\inf_{0\le w\le u}P(w)}.
\end{equation}
It follows that the sequence
\begin{equation*}
\lambda_{v}'=\frac{\lambda_{v}}{\lambda_{1}}
\end{equation*}
is bounded, and normalized so that $\lambda_{1}'=1$.\\

\begin{proof}[of Lemma \ref{L:weightedsieve}]
Let us start by showing that the number of elements of $\mathcal{A}$ that are divisible by the square of a prime $p$ with $y\le p<z$ will be relatively small. Suppose that $y=x^{\frac{1}{\alpha}}$, and $z=x^{\frac{1}{U}}$. Nagel \cite{Nag} has shown that $\rho(p^2)\le \kappa \Delta_{L}^2$. Therefore,
\begin{align*}
\sum_{y\le p<z}\left|\mathcal{A}_{p^2}\right|\ll x\sum_{p\ge y}\frac{1}{p^2}+\sum_{p<z}1\ll\frac{x}{y}+z\ll x^{1-\frac{1}{\alpha}},
\end{align*}
if $\frac{1}{U}<1-\frac{1}{\alpha}$, a condition that will be satisfied in the present application. The remaining set will be denoted by
\begin{align*}
\mathcal{A}'=\mathcal{A}-\bigcup_{y\le p<z}\mathcal{A}_{p^2}.
\end{align*}

Let $n$ denote a generic element of $\mathcal{A}'$ with $(n,P(y))=1$. If $n$ contains a repeated prime factor $p$, say, then $p\ge z$, that is, $1-\frac{\log p}{\log z}\le 0$. It follows that
\begin{equation}\label{B:richertweightbound}
\sum_{\substack{y\le p<z\\p\mid n}}\left(1-\frac{\log p}{\log z}\right)\ge\sideset{}{^*}\sum_{\substack{p\ge y\\p\mid n}}\left(1-\frac{\log p}{\log z}\right)=\Omega(n)-\frac{U\log\left|n\right|}{\log x}\ge\Omega(n)-U(\kappa+\varepsilon),
\end{equation}
where the $\ast$ denotes summation with appropriate multiplicity. The last inequality follows from the fact that
\begin{align*}
\max_{n\in\mathcal{A}}\left|n\right|\le x^{\kappa+\varepsilon}
\end{align*}
for any $\varepsilon>0$, provided $x$ is sufficiently large.\\
\indent Let $r_{\kappa}$ be a natural number such that $r_{\kappa}+1>U(\kappa+\varepsilon)$, and choose
\begin{equation*}
b=r_{\kappa}+1-U(\kappa+\varepsilon).
\end{equation*}
Using this choice of $b$ together with the bound in \eqref{B:richertweightbound}, we have that
\begin{align*}
\sum_{\substack{n\in\mathcal{A}'\\(n,P(y))=1}}\left(\sum_{\substack{d\mid n\\d\mid P(z)}}a_{d}\right)\left(\sum_{\substack{\nu\mid n\\\nu\mid P(z')}}\lambda_{\nu}'\right)^{2}&\le\sum_{\substack{n\in\mathcal{A}'\\(n,P(y))=1}}\left(b-\Omega(n)+U(\kappa+\varepsilon)\right)\left(\sum_{\substack{\nu\mid n\\\nu\mid P(z')}}\lambda_{\nu}'\right)^{2}\\
&=\sum_{\substack{n\in\mathcal{A}'\\(n,P(y))=1}}\left(r_{\kappa}+1-\Omega(n)\right)\left(\sum_{\substack{\nu\mid n\\\nu\mid P(z')}}\lambda_{\nu}'\right)^{2}.
\end{align*}
If $\Omega(n)\ge r_{\kappa}+1$, then the weight for $n$ in the last sum is either negative or zero. Thus
\begin{align*}
\sum_{\substack{n\in\mathcal{A}'\\(n,P(y))=1}}\left(\sum_{\substack{d\mid n\\d\mid P(z)}}a_{d}\right)\left(\sum_{\substack{\nu\mid n\\\nu\mid P(z')}}\lambda_{\nu}'\right)^{2}&\le\sum_{\substack{n\in\mathcal{A}', \Omega(n)\le r_{\kappa}\\(n,P(y))=1}}\left(r_{\kappa}+1-\Omega(n)\right)\left(\sum_{\substack{\nu\mid n\\\nu\mid P(z')}}\lambda_{\nu}'\right)^{2}\\
&\le \sum_{\substack{n\in\mathcal{A}',\Omega(n)\le r_{\kappa}\\(n,P(y))=1}}(r_{\kappa}+1)\left(\sum_{\substack{\nu\mid n\\\nu\mid P(z')}}\lambda_{\nu}'\right)^{2}.
\end{align*}
Observe that if $\Omega(n)\le r_{\kappa}$, then \eqref{B:lambdabound} implies that
\begin{align*}
\left(\sum_{\substack{\nu\mid n\\\nu\mid P(z')}}\lambda_{\nu}'\right)^{2}\ll\left(\sum_{\nu\mid n}\mu^{2}(\nu)\right)^{2}\ll 4^{r_{\kappa}}.
\end{align*}
Also, note that if $(n,P(y))=1$, then
\begin{align*}
\sum_{\substack{d\mid n\\d\mid P(z)}}a_{d}\ll r_{\kappa},\\
\end{align*}
and so
\begin{align*}
\sum_{\substack{n\in\mathcal{A}\\\Omega(n)\le r_{\kappa}}}1&\ge\sum_{\substack{n\in\mathcal{A}',\Omega(n)\le r_{\kappa}\\(n,P(y))=1}}1\gg4^{-r_{\kappa}}\sum_{\substack{n\in\mathcal{A}',\Omega(n)\le r_{\kappa}\\(n,P(y))=1}}\left(\sum_{\substack{\nu\mid n\\\nu\mid P(z')}}\lambda_{\nu}'\right)^{2}\\
&\gg\frac{4^{-r_{\kappa}}}{r_{\kappa}}\sum_{\substack{n\in\mathcal{A}'\\(n,P(y))=1}}\left(\sum_{\substack{d\mid n\\d\mid P(z)}}a_{d}\right)\left(\sum_{\substack{\nu\mid n\\\nu\mid P(z')}}\lambda_{\nu}'\right)^{2}\\
&\gg\sum_{\substack{n\in\mathcal{A}\\(n,P(y))=1}}\left(\sum_{\substack{d\mid n\\d\mid P(z)}}a_{d}\right)\left(\sum_{\substack{\nu\mid n\\\nu\mid P(z')}}\lambda_{\nu}'\right)^{2}+O\left(\sum_{y\le p<z}\left|\mathcal{A}_{p^2}\right|\right).
\end{align*}
Recall that if $(n,P(y))>1$, then the sum over $a_{d}$ is non-positive, so the above is
\begin{align*}
&\gg\sum_{n\in\mathcal{A}}\left(\sum_{\substack{d\mid n\\d\mid P(z)}}a_{d}\right)\left(\sum_{\substack{\nu\mid n\\\nu\mid P(z')}}\lambda_{\nu}'\right)^{2}+O\left(x^{1-\frac{1}{\alpha}}\right).
\end{align*}
In summary, we have that
\begin{align*}
\left|\left\{n\in\mathcal{A}:\Omega(n)\le r_{\kappa}\right\}\right|&\gg\frac{1}{\lambda_{1}^2}\left(x \mathfrak{S}_{\mathcal{A}}+\mathfrak{E}_{\mathcal{A}}\right)+O\left(x^{1-\frac{1}{\alpha}}\right).
\end{align*}
Using the bounds from \eqref{B:firstlambdabound} and \eqref{B:Vbound} in the above inequality leads to
\begin{align*}
\left|\left\{n\in\mathcal{A}:\Omega(n)\le r_{\kappa}\right\}\right|&\gg \frac{x}{\log^{\kappa}x}\Big(V(z')\mathfrak{S}_{\mathcal{A}}+O\left(\frac{V(z')\mathfrak{E}_{\mathcal{A}}}{x}\right)+o(1)\Big).
\end{align*}
Now we dispose of the error term involving $\mathfrak{E}_{\mathcal{A}}$.\\
\indent Let $\nu(m)$ denote the number of distinct prime divisors of $m$. Since $a_{d}=0$ for $d\ge z$, and $\lambda_{\nu}=0$ for $\nu\ge\xi$, we have
\begin{equation*}
\sum_{\substack{d\mid P(z)\\\nu_{1},\nu_{2}\mid P(z')}}a_{d}\lambda_{\nu_{1}}\lambda_{\nu_{2}}\mathcal{R}_{\left[d,\nu_{1},\nu_{2}\right]}\ll\sum_{\substack{m<z\xi^{2}\\m\mid P(z)}}\left|\mathcal{R}_{m}\right|\sum_{\substack{d,\nu_{1},\nu_{2}\\\left[d,\nu_{1},\nu_{2}\right]=m}}1=\sum_{\substack{m<z\xi^{2}\\m\mid P(z)}}7^{\nu(m)}\left|\mathcal{R}_{m}\right|.
\end{equation*}
Thus,
\begin{align*}
\mathfrak{E}_{\mathcal{A}}&\le z\xi^{2}\sum_{m\mid P(z)}\frac{7^{\nu(m)}\rho(m)}{m}\ll z\xi^{2}\prod_{p<z}\left(1+\frac{7\rho(p)}{p}\right)\\
&\ll z\xi^{2}\prod_{p<z}\left(1+\frac{\rho(p)}{p}\right)^{7}\ll z\xi^{2}\prod_{p<z}\left(1-\frac{\rho(p)}{p}\right)^{-7}\ll\frac{z\xi^{2}}{V(z)^{7}}.
\end{align*}
Choosing $z\xi^2=x^{1-\delta}$, for any $\delta>0$, we have
\begin{align*}
\left|\left\{n\in\mathcal{A}:\Omega(n)\le r_{\kappa}\right\}\right|&\gg \frac{x}{\log^{\kappa}x}\Big(V(z')\mathfrak{S}_{\mathcal{A}}+o(1)\Big),
\end{align*}
which was the desired result.
\end{proof}
\section{An innovation of Halberstam and Richert}\label{innovation}
Halberstam and Richert \cite{HR74} considered expressions of the form
\begin{equation*}
\sum_{n\in\mathcal{A}}\left(\sum_{\substack{d|n\\d\mid P(z)}}a_d\right)\left(\sum_{\substack{\nu|n\\\nu\mid P(z')}}\lambda_\nu \right)^2,
\end{equation*}
where $z$ and $z'$ are possibly distinct. If we use the identity in \eqref{E:Sieve} to decompose this expression, then we have
\begin{equation}\label{E:rawexpression}
\mathfrak{S}_{\mathcal{A}}>\sum_{\substack{m<\xi\\m\mid P(z')}}\sum_{d<z}\frac{\mu^{2}(m)}{f'(m)}\frac{a_{d}}{f(d)}\left(\sum_{r\mid d}\mu(r)\zeta_{rm}\right)^{2},
\end{equation}
upon discarding the condition that $(d,m)=1$, since discarding this condition introduces negative contributions to the sum in light of Richert's weights in \eqref{W:richertweights}. The expression on the right-hand side of \eqref{E:rawexpression} can be rewritten using Riemann-Stieltjes integration. The integrators that we will use are
\begin{equation*}
G(r,z')=\sum_{\substack{m<\xi\\m\mid P(z')}}\frac{\mu^2(m)}{f'(m)},
\end{equation*}
and
\begin{equation*}
H(s)=\sum_{p<s}\frac{\log p}{f(p)}.
\end{equation*}
Thus, inequality \eqref{E:rawexpression} is equivalent to
\begin{equation}\label{E:rawexpression2}
\mathfrak{S}_{\mathcal{A}}>b\mathcal{S}_{1}-b\mathcal{S}_{2}-\mathcal{S}_{3},
\end{equation}
where
\begin{align*}
\mathcal{S}_{1}&=\int_{1^{-}}^{\xi}\zeta_{r}^{2}dG(r,z'),\\
\mathcal{S}_{2}&=\int_{1^{-}}^{\xi}\int_{1^{-}}^{y}\left(\zeta_{r}-\zeta_{rs}\right)^{2}\frac{dH(s)}{\log s}dG(r,z'),\\
\mathcal{S}_{3}&=\int_{1^{-}}^{\xi}\int_{y^{-}}^{z}\left(\zeta_{r}-\zeta_{rs}\right)^{2}\left(1-\frac{\log s}{\log z}\right)\frac{dH(s)}{\log s}dG(r,z').
\end{align*}
We pause here to record the asymptotic formulas
\begin{equation*}
G(r,z')\sim\frac{j_{\kappa}\left(\frac{\log r}{\log z'}\right)}{V(z')},
\end{equation*}
and
\begin{equation*}
H(s)\sim\kappa\log s.
\end{equation*}
The function $j_{\kappa}(w)$ is the continuous solution of the differential delay equation
\begin{equation*}
w j_{\kappa}'(w)=\kappa j_{\kappa}(w)-\kappa j_{\kappa}(w-1),
\end{equation*}
and is defined for $0<w\le 1$ by
\begin{equation*}
j_{\kappa}(w)=c_{\kappa}w^{\kappa},
\end{equation*}
where $c_{\kappa}=\frac{e^{-\gamma\kappa}}{\Gamma(\kappa+1)}$, and $j_{\kappa}(w)=0$ if $w\le0$. More specifically, if one regards $\kappa$ and $u:=\frac{\log \xi}{\log z'}$ as fixed, then one has
\begin{equation}\label{E:Gsum}
G(r,z')=\frac{j_{\kappa}\left(\frac{\log r}{\log z'}\right)}{V(z')}\left(1+O\left(\frac{1}{\log z'}\right)\right),
\end{equation}
and
\begin{equation}\label{E:HSum2}
H(s)=\kappa\log s+O\left(1\right).
\end{equation}
The formula in \eqref{E:HSum2} is merely our assumed density hypothesis of $\rho(p)$ in \eqref{E:HSum}. On the other hand, the bound in \eqref{E:Gsum} is a consequence of
\begin{lemma}\label{T:HRlemma}
For any $\tau=\frac{\log r}{\log z}>0$, we have
\[
\frac{1}{G(r,z)}=V(z)\left(\frac{1}{j_{\kappa}(\tau)}+O\left(\frac{\tau^{2\kappa+1}}{\log z}\right)\right).
\]
\end{lemma}
Lemma \ref{T:HRlemma} is discussed in some detail in Halberstam and Richert \cite[see Section 4 on p.197]{HR74}. Now, in view of the asymptotic formulas above in \eqref{E:HSum2} and \eqref{E:Gsum}, we have
\begin{align*}
\mathcal{S}_{1}&\sim \frac{1}{V(z')}\int_{1}^{\xi}\zeta_{r}^{2}dj\left(\frac{\log r}{\log z'}\right),\\
\mathcal{S}_{2}&\sim \frac{\kappa}{V(z')}\int_{1}^{\xi}\int_{1}^{y}\left(\zeta_{r}-\zeta_{rs}\right)^{2}\frac{d\log s}{\log s}dj\left(\frac{\log r}{\log z'}\right),\\
\mathcal{S}_{3}&\sim \frac{\kappa}{V(z')}\int_{1}^{\xi}\int_{y}^{z}\left(\zeta_{r}-\zeta_{rs}\right)^{2}\left(1-\frac{\log s}{\log z}\right)\frac{d\log s}{\log s}dj\left(\frac{\log r}{\log z'}\right),
\end{align*}
where the error term in each of these asymptotic relations is of order at most $\left(V(z')\log z'\right)^{-1}$. Putting these asymptotic formulas together with \eqref{E:rawexpression2}, we have that
\begin{equation*}
\mathfrak{S}_{\mathcal{A}}\gtrsim\frac{1}{V(z')}\left(b\mathcal{I}^{*}_{1}-\kappa\mathcal{I}^{*}_{2}+\kappa\mathcal{I}^{*}_{3}-b\kappa\mathcal{I}^{*}_{4}\right),
\end{equation*}
where
\begin{align*}
\mathcal{I}^{*}_{1}&=\int_{1}^{\xi}\zeta_{r}^2dj\left(\frac{\log r}{\log z'}\right),\\
\mathcal{I}^{*}_{2}&=\int_{1}^{\xi}\int_{1}^{z}\left(\zeta_{r}-\zeta_{rs}\right)^{2}\left(1-\frac{\log s}{\log z}\right)\frac{d\log s}{\log s}dj\left(\frac{\log r}{\log z'}\right),\\
\mathcal{I}^{*}_{3}&=\int_{1}^{\xi}\int_{1}^{y}\left(\zeta_{r}-\zeta_{rs}\right)^{2}\left(1-\frac{\log s}{\log z}\right)\frac{d\log s}{\log s}dj\left(\frac{\log r}{\log z'}\right),\\
\mathcal{I}^{*}_{4}&=\int_{1}^{\xi}\int_{1}^{y}\left(\zeta_{r}-\zeta_{rs}\right)^{2}\frac{d\log s}{\log s}dj\left(\frac{\log r}{\log z'}\right).
\end{align*}
Let us suppose that $z=x^{1/U}$, $z'=x^{1/V}$, $u=\log \xi/\log z'$. Recall from equation \eqref{E:polyP} that
\begin{equation*}
\zeta_{r}=P\left(\frac{\log\xi/r}{\log z'}\right)
\end{equation*}
when $r$ is squarefree, $r<\xi$, and $r\mid P(z')$. At this point, it is convenient to define
\begin{equation*}
P^{*}(w)=
\begin{cases}
P(w)& \text{if $w\ge0$,}\\
0& \text{otherwise.}
\end{cases}
\end{equation*}
Making the change of variables $v=\log r/\log z'$ and $t=\log s/\log z'$, the integrals above can be rewritten as
\begin{align*}
\mathcal{I}^{*}_{1}&=\int_{0}^{u}P^{*}\left(u-v\right)^{2}j'(v)dv,\\
\mathcal{I}^{*}_{2}&=\int_{0}^{u}\int_{0}^{V/U}\left(P^{*}\left(u-v\right)-P^{*}\left(u-v-t\right)\right)^{2}\left(1-t\frac{U}{V}\right)\frac{dt}{t}j'(v)dv,\\
\mathcal{I}^{*}_{3}&=\int_{0}^{u}\int_{0}^{V/\alpha}\left(P^{*}\left(u-v\right)-P^{*}\left(u-v-t\right)\right)^{2}\left(1-t\frac{U}{V}\right)\frac{dt}{t}j'(v)dv\\
\mathcal{I}^{*}_{4}&=\int_{0}^{u}\int_{0}^{V/\alpha}\left(P^{*}\left(u-v\right)-P^{*}\left(u-v-t\right)\right)^{2}\frac{dt}{t}j'(v)dv.
\end{align*}
Next, let $w=u-v$, and $l=V/U\ge1$, so that
\begin{align*}
\mathcal{I}^{*}_{1}&=\int_{0}^{u}P^{*}\left(w\right)^{2}j'(u-w)dw,\\
\mathcal{I}^{*}_{2}&=\int_{0}^{u}\int_{0}^{l}\left(P^{*}\left(w\right)-P^{*}\left(w-t\right)\right)^{2}\left(1-\frac{t}{l}\right)\frac{dt}{t}j'(u-w)dw,\\
\mathcal{I}^{*}_{3}&=\int_{0}^{u}\int_{0}^{V/\alpha}\left(P^{*}\left(w\right)-P^{*}\left(w-t\right)\right)^{2}\left(1-\frac{t}{l}\right)\frac{dt}{t}j'(u-w)dw\\
\mathcal{I}^{*}_{4}&=\int_{0}^{u}\int_{0}^{V/\alpha}\left(P^{*}\left(w\right)-P^{*}\left(w-t\right)\right)^{2}\frac{dt}{t}j'(u-w)dw.
\end{align*}
These last two integrals can be made as small as we like provided we take $\alpha$ sufficiently large, and so we have
\begin{align*}
\mathfrak{S}_{\mathcal{A}}&\gtrsim\frac{1}{V(z')}\left(b\mathcal{I}^{*}_{1}-\kappa\mathcal{I}^{*}_{2}\right).
\end{align*}
Now, let us assume that $u\le l$. To account for the fact that $P^{*}(w-t)=0$ when $w\le t$, we split the range of the innermost integral appearing in $\mathcal{I}^{*}_{2}$. This proves
\begin{lemma}\label{L:maintermlemma}
Suppose that $\mathcal{A}=\left\{L(n):n\le x\right\}$ and that $L(n)$ satisfies the hypotheses of Theorem \ref{Th:1}. Let $y=x^{\frac{1}{\alpha}}$, $z=x^{\frac{1}{U}}$, $z'=x^{\frac{1}{V}}$, $\xi^{\frac{1}{u}}=z'$, and $l=\frac{V}{U}\ge1$. Then, for all sufficiently large $\alpha$ and $x$,
\begin{equation*}
\mathfrak{S}_{\mathcal{A}}\gtrsim\frac{1}{V(z')}\left(b\mathcal{I}_{1}-\kappa\mathcal{I}_{2}-\kappa\mathcal{I}_{3}\right),
\end{equation*}
where
\begin{align*}
\mathcal{I}_{1}&=\int_{0}^{u}P\left(w\right)^{2}j'(u-w)dw,\\
\mathcal{I}_{2}&=\int_{0}^{u}\int_{0}^{w}\left(P\left(w\right)-P\left(w-t\right)\right)^{2}\left(1-\frac{t}{l}\right)\frac{dt}{t}j'(u-w)dw,\\
\mathcal{I}_{3}&=\int_{0}^{u}\int_{w}^{l}P\left(w\right)^{2}\left(1-\frac{t}{l}\right)\frac{dt}{t}j'(u-w)dw.
\end{align*}
\end{lemma}

Recall from Lemma \ref{L:weightedsieve} that the error terms will be kept under control if
\begin{align*}
z\xi^{2}=x^{1-\delta},
\end{align*}
for any $\delta>0$, and since
\begin{align*}
z\xi^{2}=zz'^{2u}=x^{\frac{1}{U}}x^{\frac{2u}{V}}=x^{\frac{1}{U}+\frac{2u}{V}},
\end{align*}
we choose
\begin{align*}
\frac{1}{U}+\frac{2u}{V}=1-\delta,
\end{align*}
or equivalently,
\begin{equation}\label{E:choiceforU}
U=1+\frac{2u}{l}+O(\delta).
\end{equation}
Following Richert, we choose
\begin{equation}\label{choiceofb}
b=r_{\kappa}+1-\left(\kappa+\varepsilon\right)U=r_{\kappa}+1-\kappa\left(1+\frac{2u}{l}\right)+O(\varepsilon'),
\end{equation}
where $\varepsilon'>0$ can be made arbitrarily small.

\section{An application of the saddle point method}

To obtain an improvement in the bound for $r_{\kappa}$ when $\kappa$ is large, one must handle integrals of the form
\begin{equation}\label{E:moment1}
\mathcal{J}_{1}(i)=\int_{0}^{u}w^{i}j'\left(u-w\right)dw,
\end{equation}
and
\begin{equation}\label{E:moment2}
\mathcal{J}_{2}(i)=\int_{0}^{u}w^{i}\log w\ j'\left(u-w\right)dw,
\end{equation}
when $u$ is around $\kappa$. Selberg encountered integrals of this form when obtaining an asymptotic formula for the sifting limit of his lower bound sieve. His calculations can be found in \cite[equation (14.23)]{Selberg}. The key to evaluating such integrals is an asymptotic formula for $j'(u-w)$ obtained by the saddle point method. With $u=\kappa-1/3-d$, Selberg showed that
\begin{equation}\label{E:j' estimate}
j'(u-w)=\frac{1}{\sqrt{\pi\kappa}}e^{-\frac{w^2}{\kappa}}\left(1-2d\frac{w}{\kappa}-\frac{4}{9}\frac{w^3}{\kappa^2}+O\left(\frac{1}{\kappa}+\frac{w^6}{\kappa^4}\right)\right),
\end{equation}
for $0\le w\le \kappa^{3/5}$. This estimate can be used to prove
\begin{lemma}\label{L:selbergintegral}
Suppose $u=\kappa-1/9$. Then, we have
	\begin{equation}\label{newJ10}
\mathcal{J}_{1}(0)=\frac{1}{2}+O\left(\frac{1}{\kappa}\right),
	\end{equation}
	\begin{equation}\label{newJ11}
\mathcal{J}_{1}(1)=\frac{1}{2}\sqrt{\frac{\kappa}{\pi}}-\frac{1}{18}+O\left(\frac{1}{\sqrt{\kappa}}\right),
	\end{equation}	
and
	\begin{equation}\label{newJ2}
\mathcal{J}_{2}(0)=\frac{1}{4}\log\kappa+\frac{1}{4}\Psi\left(\frac{1}{2}\right)-\frac{1}{9\sqrt{\pi\kappa}}+O\left(\frac{\log \kappa}{\kappa} \right),
	\end{equation}
where $\Psi(z)$ is the digamma function,
	\begin{equation}\label{Psi}
\Psi(z)=\frac{\Gamma'(z)}{\Gamma(z)}.
	\end{equation}
\end{lemma}
Thus, setting $d=-2/9$, we calculate that
\begin{equation}\label{ratio1Th1}
\frac{\mathcal{J}_{1}(1)}{\mathcal{J}_{1}(0)}=\sqrt{\frac{\kappa}{\pi}}-\frac{1}{9}+O\left(\frac{1}{\sqrt{\kappa}}\right),
\end{equation}
and
	\begin{equation}\label{ratio2Th1}
\frac{\mathcal{J}_{2}(0)}{\mathcal{J}_{1}(0)}=\frac{1}{2}\log\kappa+\frac{1}{2}\Psi\left(\frac{1}{2}\right)-\frac{2}{9\sqrt{\pi\kappa}}+O\left(\frac{\log\kappa}{\kappa}\right).
	\end{equation}
If the reader wishes to use these ratios to calculate the main term of the bound for $r_{\kappa}$ in Theorem \ref{Th:1}, then the remainder of this section can be skipped.

\begin{proof}
We will prove \eqref{newJ2} and leave the proof of \eqref{newJ10} and \eqref{newJ11} to the reader, since it follows in exactly the same manner. First, for any $i\ge0$,
\begin{equation}\label{E:log gamma integral}
\int_{0}^{\infty}w^{i}\log w\ e^{-\frac{w^2}{\kappa}}dw=\frac{1}{4}\Gamma\left(\frac{i+1}{2}\right)\kappa^{\frac{i+1}{2}}\log \kappa+\frac{1}{4}\Gamma\left(\frac{i+1}{2}\right)\kappa^{\frac{i+1}{2}}\Psi\left(\frac{i+1}{2}\right).
\end{equation}
This is easily seen by performing the change of variable $t=w^2/\kappa$, for
\begin{align*}
\int_{0}^{\infty}w^{i}\log w\ e^{-\frac{w^2}{\kappa}}dw&=\frac{1}{4}\kappa^{\frac{i+1}{2}}\int_{0}^{\infty}t^{\frac{i-1}{2}}\log(\kappa t)\ e^{-t}dt\\
&=\frac{1}{4}\kappa^{\frac{i+1}{2}}\log\kappa\int_{0}^{\infty}t^{\frac{i-1}{2}}\ e^{-t}dt+\frac{1}{4}\kappa^{\frac{i+1}{2}}\int_{0}^{\infty}t^{\frac{i-1}{2}}\log t\ e^{-t}dt,
\end{align*}
and these last two integrals are the integral representations of $\Gamma\left(\frac{i+1}{2}\right)$ and $\Gamma'\left(\frac{i+1}{2}\right)$. The integral representation of $\Gamma'\left(\frac{i+1}{2}\right)$ is obtained by an application of differentiation under the integral, and we use \eqref{Psi} to write this in terms of the digamma function, $\Psi\left(\frac{i+1}{2}\right)$.\\

We will also make use of an estimate for
\begin{equation}\label{E:log gamma tail}
\int_{\kappa^{3/5}}^{\infty}w^{i}\log w\ e^{-\frac{w^2}{\kappa}}dw\ll\kappa^{(2+3i)/5}\log \kappa\ e^{-\kappa^{1/5}},
\end{equation}
which follows from an elementary estimate of the incomplete gamma function
\begin{equation*}
\Gamma\left(s,x\right)=\int_{x}^{\infty}t^{s-1}e^{-t}dt\sim x^{s-1}e^{-x},
\end{equation*}
valid as long as $s=o(x)$. A proof of this estimate is easy to supply since, upon integrating by parts, we have
\begin{equation*}
\Gamma\left(s,x\right)=x^{s-1}e^{-x}+(s-1)\int_{x}^{\infty}t^{s-2}e^{-t}dt=x^{s-1}e^{-x}+O\left(\frac{s-1}{x}\Gamma\left(s,x\right)\right).
\end{equation*}
Returning to \eqref{E:log gamma tail}, we have
\begin{equation*}
\int_{\kappa^{\frac{3}{5}}}^{\infty}w^{i}e^{-\frac{w^2}{\kappa}}\log w\ dw=\frac{1}{4}\kappa^{\frac{i+1}{2}}\Gamma\left(\frac{i+1}{2},\kappa^{\frac{1}{5}}\right)\log\kappa+\frac{1}{4}\kappa^{\frac{i+1}{2}}\int_{\kappa^{\frac{1}{5}}}^{\infty}t^{\frac{i-1}{2}}e^{-t}\log t\ dt,
\end{equation*}
and
\begin{align*}
\int_{\kappa^{\frac{1}{5}}}^{\infty}t^{\frac{i-1}{2}}e^{-t}\log t\ dt&=\int_{\kappa^{\frac{1}{5}}}^{\kappa}t^{\frac{i-1}{2}}e^{-t}\log t\ dt+\int_{\kappa}^{\infty}t^{\frac{i-1}{2}}e^{-t}\log t\ dt\\
&\ll\log \kappa \int_{\kappa^{\frac{1}{5}}}^{\kappa}t^{\frac{i-1}{2}}e^{-t} dt+\int_{\kappa}^{\infty}t^{\frac{i-1}{2}+\varepsilon}e^{-t}dt\\
&\ll\Gamma\left(\frac{i+1}{2},\kappa^{\frac{1}{5}}\right)\log \kappa+\Gamma\left(\frac{i+1}{2}+\varepsilon,\kappa\right)\\
&\ll\Gamma\left(\frac{i+1}{2},\kappa^{\frac{1}{5}}\right)\log \kappa.
\end{align*}
With the estimates above, we can begin analysis of $\mathcal{J}_{2}(0)$. We first split the range of integration to obtain
\begin{equation}\label{J2}
\mathcal{J}_{2}(0)=\mathcal{J}_{2,1}(0)+\mathcal{J}_{2,2}(0),
\end{equation}
say, where
\begin{equation*}
\mathcal{J}_{2,1}(0)=\int_{0}^{\kappa^{3/5}}\log w\ j'\left(u-w\right)dw,
\end{equation*}
and
\begin{equation*}
\mathcal{J}_{2,2}(0)=\int_{\kappa^{3/5}}^{u}\log w\ j'\left(u-w\right)dw.
\end{equation*}

Next, we dispose of $\mathcal{J}_{2,2}(0)$ using integration by parts together with the inequality
\begin{equation*}
j(u-w)\le e^{-\frac{w^2}{\kappa}},
\end{equation*}
valid for $\kappa^{3/5}<w\le u$ \cite{Selberg}. This shows that
\begin{equation}\label{J22}
\mathcal{J}_{2,2}(0)=-\int_{\kappa^{3/5}}^{u}\log w\ dj(u-w)\ll\log \kappa\ e^{-\kappa^{1/5}}.
\end{equation}

Moving on, we plug the asymptotic formula for $j'(u-w)$ given in \eqref{E:j' estimate} into $\mathcal{J}_{2,1}(0)$, and distinguish
\begin{equation}\label{J21}
\mathcal{J}_{2,1}(0)=\mathcal{J}_{2,1,1}(0)+\mathcal{J}_{2,1,2}(0),
\end{equation}
where
\begin{equation*}
\mathcal{J}_{2,1,1}(0)=\int_{0}^{\kappa^{3/5}}M\left(w,\kappa\right)\log w\ e^{-\frac{w^2}{\kappa}}dw,
\end{equation*}
\begin{equation*}
\mathcal{J}_{2,1,2}(0)=\int_{0}^{\kappa^{3/5}}E\left(w,\kappa\right)\log w\ e^{-\frac{w^2}{\kappa}}dw,
\end{equation*}
\begin{equation*}
M(w,\kappa)=\frac{1}{\sqrt{\pi\kappa}}\left(1+\frac{4}{9}\frac{w}{\kappa}-\frac{4}{9}\frac{w^3}{\kappa^2}\right),
\end{equation*}
and
\begin{equation*}
E(w,\kappa)\ll\frac{1}{\sqrt{\pi\kappa}}\left(\frac{1}{\kappa}+\frac{w^6}{\kappa^4}\right).
\end{equation*}

Observe that, using \eqref{E:log gamma integral}, we have
\begin{equation}\label{J212}
\mathcal{J}_{2,1,2}(0)\ll\int_{0}^{\infty}E(w,\kappa)\log w\ e^{-\frac{w^2}{\kappa}}dw\ll\frac{\log\kappa}{\kappa},
\end{equation}
and
\begin{align*}
\mathcal{J}_{2,1,1}(0)&=\int_{0}^{\infty}M(w,\kappa)\log w\ e^{-\frac{w^2}{\kappa}}dw+O\left(\int_{\kappa^{3/5}}^{\infty}M(w,\kappa)\log w\ e^{-\frac{w^2}{\kappa}}dw\right)\\
&=\frac{1}{4}\log\kappa+\frac{1}{4}\Psi\left(\frac{1}{2}\right)-\frac{1}{9\sqrt{\pi\kappa}}+O\left(\int_{\kappa^{3/5}}^{\infty}M(w,\kappa)\log w\ e^{-\frac{w^2}{\kappa}}dw\right).
\end{align*}
Finally, using \eqref{E:log gamma tail},
\begin{equation}\label{J211}
\mathcal{J}_{2,1,1}(0)=\frac{1}{4}\log\kappa+\frac{1}{4}\Psi\left(\frac{1}{2}\right)-\frac{1}{9\sqrt{\pi\kappa}}+O\left(\kappa^{-1/10}\log \kappa\ e^{-\kappa^{1/5}}\right).
\end{equation}
The asymptotic formula for $\mathcal{J}_{2}(0)$ in \eqref{newJ2} follows by combining \eqref{J2}, \eqref{J22}, \eqref{J21}, \eqref{J212}, and \eqref{J211}.
\end{proof}
\section{Proof of Theorem 1}
In this section, we prove Theorem \ref{Th:1}. We will choose $u=\kappa-1/9$, and $P(w)=1$. Thus, the improvement over other authors is attributed to the large choice of $u$, for which we use the asymptotic formula for $j'(u-w)$ in \eqref{E:j' estimate}. The device of Halberstam and Richert allows us to choose this large $u$ and still keep $U$ small by taking $z'$ smaller than $z$.
\begin{proof}[of Theorem \ref{Th:1}]
Choose $u=\kappa-\frac{1}{9}$, and $P(w)=1$, and observe that Theorem \ref{Th:1} follows from Lemma \ref{L:weightedsieve} and Lemma \ref{L:maintermlemma} if
\begin{align*}
b\int_{0}^{u}j'(u-w)dw-\kappa\int_{0}^{u}\int_{w}^{l}\left(1-\frac{t}{l}\right)\frac{dt}{t}j'(u-w)dw>0.
\end{align*}
Computing the innermost integral, this inequality becomes
\begin{equation*}
b>\kappa(\log l-1)-\kappa\frac{\displaystyle\int_{0}^{u}\log w j'(u-w)dw}{\displaystyle\int_{0}^{u}j'(u-w)dw}+\frac{\kappa}{l}\frac{\displaystyle\int_{0}^{u}wj'(u-w)dw}{\displaystyle\int_{0}^{u}j'(u-w)dw}.
\end{equation*}
Using formulas \eqref{ratio1Th1} and \eqref{ratio2Th1} for these ratios of integrals, we see that
\begin{equation*}
\frac{\displaystyle\int_{0}^{u}\log w\  j'(u-w)dw}{\displaystyle\int_{0}^{u}j'(u-w)dw}=\frac{1}{2}\log\kappa+\frac{1}{2}\Psi\left(\frac{1}{2}\right)-\frac{2}{9\sqrt{\pi\kappa}}+O\left(\frac{\log\kappa}{\kappa}\right),
\end{equation*}
and
\begin{equation*}
\frac{\displaystyle\int_{0}^{u}wj'(u-w)dw}{\displaystyle\int_{0}^{u}j'(u-w)dw}=\sqrt{\frac{\kappa}{\pi}}-\frac{1}{9}+O\left(\frac{1}{\sqrt{\kappa}}\right).
\end{equation*}
Now, plugging these ratios into the inequality above and recalling that
\begin{equation*}
b=r+1-\kappa\left(1+\frac{2u}{l}\right)+O\left(\varepsilon'\right),
\end{equation*}
we have
\begin{equation*}
r_{\kappa}>\frac{1}{2}\kappa\log\kappa+\left(\frac{2\kappa}{l}-\frac{1}{2}\Psi\left(\frac{1}{2}\right)+\log\frac{l}{\kappa}\right)\kappa+\left(\frac{2}{9}+\frac{\kappa}{l}\right)\sqrt{\frac{\kappa}{\pi}}+O\left(\frac{\kappa}{l}+\frac{\log\kappa}{\sqrt{\kappa}}\right).
\end{equation*}
Setting $l=2\kappa$, and using $\Psi\left(\frac{1}{2}\right)=-\gamma-2\log 2$, this becomes
\begin{equation*}
r_{\kappa}>\frac{1}{2}\kappa\log\kappa+\left(1+\frac{\gamma}{2}+\log 4\right)\kappa+\frac{13}{18}\sqrt{\frac{\kappa}{\pi}}+O\left(\log\kappa\right).
\end{equation*}
Since any $r_{\kappa}$ satisfying the inequality above must also satisfy
\begin{equation*}
r_{\kappa}>\kappa\left(1+\frac{2u}{l}\right)-1+O(\varepsilon')=2\kappa-\frac{10}{9}+O(\varepsilon'),
\end{equation*}
the proof of Theorem \ref{Th:1} is complete.
\end{proof}

\section{Acknowledgements}
I would like to thank Sid Graham for all of his help and guidance throughout this research. He continues to be a great mentor. I would also like to thank Olivier Ramar\'e for carefully reading through this manuscript, taking care to provide thoughtful comments, and suggestions for future research.

%\affiliationone{% in this example, two authors share an institution
%   C. S. Franze\\
%   Central Michigan University\\
%   Department of Mathematics\\
%   Mount Pleasant, MI  48859\\
%   U.S.A.
%   \email{franz1cs@cmich.edu}}
% Important: Do not put any empty line here.
%\affiliationtwo{% in this example, one author has two addresses}
%   T. Hird\\
%   Previous postal address where
%     the research was performed and\\
%   Country
%   \email{hird@university.ac.uk}}
% Important: Do not put any empty line here.
% Use \affiliationthree{} for any address positioned under \affiliationone
% Use \affiliationfour{}  for any address positioned under \affiliationtwo
%\affiliationthree{~} %inserts a space to make this field empty
%\affiliationfour{%
%   Current address:\\
%   Present long-term address\\
%   Country
%   \email{t.hird@institution.edu}}
%
\end{document}